\RequirePackage{fix-cm}
\documentclass[smallextended]{svjour3}       % onecolumn (second format)
\smartqed  % flush right qed marks, e.g. at end of proof
\usepackage{graphicx}
\newtheorem{result}{Approximation}

\newtheorem{prop}{Proposition}
\usepackage{algorithmic,algorithm2e,float}
\usepackage{amsmath,amssymb,euscript,amsfonts,mathrsfs,amscd, amstext, mathbbol}
\makeatletter
\renewcommand{\@algocf@capt@plain}{above}% formerly {bottom}
\makeatother
\usepackage{color}
\usepackage{mathtools}
\let\oldemptyset\emptyset
\let\emptyset\varnothing
 \usepackage{multirow}

\providecommand{\abs}[1]{\lvert#1\rvert}

\begin{document}

\title{A review of the deterministic and diffusion approximations for stochastic chemical reaction networks \thanks{P. Mozgunov and T. Jaki have received funding from the European Union`s Horizon 2020 research and innovation programme under the Marie Sklodowska-Curie grant agreement No 633567.}
}
\subtitle{ }

\titlerunning{A review of approximations for stochastic chemical reaction networks}        % if too long for running head

\author{Pavel Mozgunov        \and
	Marco Beccuti \and
	Andras Horvath \and
        Thomas Jaki \and
        Roberta Sirovich \and
        Enrico Bibbona %etc.
}

\authorrunning{P. Mozgunov et.al} % if too long for running head

\institute{Pavel Mozgunov, Thomas Jaki \at
              Medical and Pharmaceutical Statistics Research Unit, Department of Mathematics and Statistics, Lancaster University  \\
              Tel.: +447493182881\\
              \email{p.mozgunov@lancaster.ac.uk, t.jaki@lancaster.ac.uk}           %  \\
%             \emph{Present address:} of F. Author  %  if needed
           \and
           Marco Beccuti, Andras Horvath \at
           Dipartimento di Informatica\\
           Universit\`a di Torino\\
           \email{beccuti@di.unito.it}\\
           \email{horvath@di.unito.it}
           \and
           Roberta Sirovich \at
           Dipartimento di Matematica G. Peano\\
           Universit\`a di Torino\\
           \email{roberta.sirovich@unito.it}
           \and
           Enrico Bibbona \at
             DISMA, Politecnico di Torino\\
             \email{enrico.bibbona@polito.it}
}

\date{Received: date / Accepted: date}
% The correct dates will be entered by the editor

\maketitle

\begin{abstract}
This work reviews deterministic and diffusion approximations of the stochastic chemical reaction networks and explains their applications. We discuss the added value the diffusion approximation provides for systems with different phenomena, such as a deficiency and a bistability. It is advocated that the diffusion approximation can be considered as an alternative theoretical approach to study the reaction networks rather than a simulation shortcut.  We discuss two examples in which the diffusion approximation is able to catch qualitative properties of reaction networks that the deterministic model misses. We provide an explicit construction of the original process and the diffusion approximation such that the distance between their trajectories is controlled and demonstrate this construction for the examples. We also discuss the limitations and potential directions of the developments.
% \PACS{PACS code1 \and PACS code2 \and more}
% \subclass{MSC code1 \and MSC code2 \and more}
\end{abstract}

\section{Introduction}

A mathematical modelling of chemical kinetics was initiated at the beginning of the previous century. This topic has attracted an extensive attention and works of many excellent scientists formed the \emph{Reaction Network Theory} at the end of the 1980s. In this formalization, the concentration of species in the network of chemical reactions obeys deterministic laws which are encoded into systems of the non-linear ordinary differential equations (ODEs). These equations provided a rich collection of complex examples and helped to improve the theory of dynamical systems \cite{erditoth,feinberg}.

Despite the fact that the deterministic models have been sufficient for the majority of applications available at that time, it was already known that a microscopic description of chemical kinetics should have included randomness (see \cite{lente} for historical remarks).
The most popular way to describe the stochastic models of reaction networks is in terms of Continuous Time Markov Chain (CTMC). The reactions are considered as happening at random events that modify the state of the network according to the stoichiometric equations. For some time these two descriptions have been developed in parallel and using different tools: the deterministic models were investigated in the theoretical and mathematical aspects, while the stochastic models were mainly studied from a computational point, e.g. in the search of suitable simulation algorithms.  

The relation between the deterministic model and the stochastic counterparts was clarified in the works by T. Kurtz \cite{kurtz1ode,KurtzChemPhys}. It was proved that the stochastic models converge to the deterministic ones if the initial amount of molecules is large. 
This was an important theoretical breakthrough in both Chemistry and Mathematics since it is appeared that deterministic and stochastic models are not independent alternative modelling frameworks. In fact, the deterministic model is an approximation of the stochastic one which is of the key importance when the system is large. Indeed, one of the practical problems of stochastic modelling is that for a system with a large amount of molecules, reactions can be so frequent that even a numerical simulation becomes computationally infeasible. Thus, the value of approximations which are easier to handle either numerically or theoretically cannot be underestimated \cite{approxreview}.

On the other hand, the deterministic approximation can lose an important information in many stochastic systems. While it usually provides a good approximation of the process mean value, it ignores completely other properties, for instance, variance, bimodality, tail behaviour, etc. T. Kurtz \cite{Kurtz1976} provided a second approximation which retains the stochastic nature by means of the diffusion process. The same equations had became popular in Chemistry under the name of \textit{Langevin equations} due to the contribution by Gillespie \cite{gillespie}. These equations have been used in many works as a computational trick to speed up simulations of the original process. While this \emph{computational} approach to the diffusion approximation proved to be fruitful, such interpretation hides in part the richness and the importance of the result by T. Kurtz \cite{Kurtz1976}. Moreover, stochastic reaction networks attracted a renewed interest recently, see e.g., \cite{AKbook,lente,santillan,ullah}. New motivations come both from the application in the system biology, demonstrating the emergence of the \emph{stochastic effects} at small scales, and from new theoretical investigations that allowed to extend mathematical results, previously known in the deterministic setting only, to the stochastic world \cite{anderson2017non,def0theorem,danielecarsten}.

The goal of this communication is to review both deterministic and diffusion approximations of the CTMC  and to explain their implications and the added value the diffusion approximation can provide for systems of the intermediate size. We emphasize that the results by T. Kurtz \cite{Kurtz1976} are constructive. It allows to give an explicit construction of the CTMC and the diffusion approximation coupled trajectories such that the uniform distance between them is controlled. To our knowledge, this fact has been never highlighted in the applied literature, while deserving to be understood better. We provide two examples in which the diffusion approximation is able to catch qualitative properties of the reaction networks that the deterministic model misses. We advocate that in the context of growing interest to the stochastic models, the diffusion approximation (or other new approximations of the same nature) has  an important role in the development of the theory and deserves to be extended for new challenges opened by the applications.

\section{Stochastic models of reaction networks and their approximations}
\subsection{Reaction networks and deficiency \label{prima}}

A reaction network is a triple $\{\mathcal{S},\mathcal{C}, \mathcal{R}\}$ such that
\begin{enumerate}
	\item $\mathcal{S}=\{S_1,\cdots, S_d\}$ is the set of species of cardinality $d$ where $d$ is finite.
	\item $\mathcal{C}$ is the set of complexes, consisting of some nonnegative {integer} linear combination of the species.
	\item $\mathcal{R}$ is a finite set of ordered couples of complexes {which is defined by the stoichiometric equations (\ref{reactionk})}.
	\end{enumerate}
A reaction network of $K$ chemical reactions is specified by stoichiometric equations 
\begin{equation}
\sum_{i=1}^d c_{ki} S_i \rightarrow \sum_{i=1}^d {c}_{ki}^\prime S_i, \ k=1,\ldots,K
\label{reactionk}
\end{equation}
meaning that the reaction consumes $\sum c_{ki} S_i$ to produce $\sum c_{ki}^\prime S_i$ where $c_{ki},{c_{ki}^\prime}$ are nonnegative  {integers}.
The definition above implies the unique directed graph {if the set of nodes coincides with the set of
complexes.} The inference of qualitative properties of the reaction network model
is based on the algebraic properties of this graph, see e.g. \cite{elisendacarsten}. We define $l_k=c^{\prime}_k-c_k$ as the reaction vector of the network where $c_k=[c_{k1},\ldots,c_{kd}]^{\rm T}$ and $c^{\prime}_k=[c^{\prime}_{k1},\ldots,c^{\prime}_{kd}]^{\rm T}$. They can be collected as the columns of the $d \times K$ \textit{stoichiometric matrix}.

One of the most important algebraic properties of the reaction network graph is the \emph{deficiency}.
Let $L$ be the number of connected components (also known as linkage classes) of the reaction graph. The subspace of $\mathbb{Z}^K$ given by
\[S=\text{span}_k\{l_k\},\]
is the \emph{stoichiometric subspace} (with dimension $\dim S$) of the network. The number of complexes in the reaction network is given by $|\mathcal{C}|$.
The \emph{deficiency} of the network is defined as the integer \[\theta=|\mathcal{C}|-L-\dim S.\] The property of non-deficiency ($\theta=0$) has important consequences on the dynamics of both deterministic and stochastic models (see Section \ref{ex} and  \emph{Deficiency-Zero} theorems \cite{def0theorem} and \cite{AKbook} for further details).

\subsection{Stochastic models of reaction networks}
%[Y_1(t),\cdots,Y_d(t)]^{\rm T}

A stochastic model of the reaction network is a Markov chain $Y(t)$ whose state space is subset of $\mathbb{N}^d$. The state vector  $s=(s_1,\cdots, s_d)$ corresponds to the number of molecules of each species available in the system. If $s_i\geq c_{ki}$ for all $i\in\{1,\cdots, d\}$, the $k^{th}$ reaction of (\ref{reactionk}) can occur, updating the network from state $s$ to  state $s + l_k$. {The occurrences of the reactions determine the jumps of the Markov chain}. The network follows the mass-action kinetics if the rate of reaction $k$ in state $s$ can be written in the form
\begin{equation}
q_{s,s+l_{k}}=\frac{\lambda_{k}}{V^{\langle c_k \rangle-1}}\prod_{i=1}^d \binom{s_i}{c_{ik}}= V\left[\frac{\lambda_{k}}{\prod_{i=1}^d c_{ik} !}\prod_{i=1}^d \left(\frac{s_i}{V}\right)^{c_{ik}} +O\left(\frac{1}{V}\right)\right]
 \label{massaction}
\end{equation}
where $\langle c_k \rangle=\sum_i c_{ik}$, $\lambda_{k}$ is a transition propensity for reaction $k$ and $V$ is the constant volume of the container in which the reactions take place.

%For monomolecular or bimolecular reaction (\ref{massaction}) gives the usual mass-action rates ~\cite{AKbook,lente}, but also includes more general models (see  Sections \ref{ex1} and \ref{ex2} for examples).

The mass-action rates \eqref{massaction} allows the Markov chain models of reaction networks to satisfy the property of the density dependence (approximately). This allows to apply the approximation results given in this communication.
	\begin{definition}
	A family of continuous time Markov Chains $\{Y^{[V]}(t)\}$ indexed by  parameter $V$ and with state spaces contained in $\mathbb{Z}^d$ is density dependent if its transition rates $q^{[V]}_{s,s+l}$ from any state $s$ to any other state $s+l$ can be written in the following form
	\begin{equation}\label{rates}
	q^{[V]}_{s,s+l}=V\; f_l \left(\frac{s}{V}\right)
	\end{equation}
	where $f_l$ is a non-negative function defined on some subset of $\mathbb{R}^d$.
	\end{definition}
\noindent Intuitively, the necessary conditions of the density dependence are (i) the linear relation of transition rates on $V$ and (ii) the dependence on the density of the population levels rather than on the population values. Then, the argument $s/V$ in (\ref{rates}) is the density associated with state $s$ and index $l$ is the vector of transitions. In case of reaction networks, the indexing parameter $V$ is the volume of the container and the process $Y^{[V]}(t)$ provides a number of molecules at time $t$. It becomes apparent from (\ref{rates}) that (\ref{massaction}) has the approximately density dependent form. In fact, it is more common to rescale the number of molecules to the concentrations. 
\begin{definition}
		For the density dependent family $\{Y^{[V]}(t)\}$ we define the family of density processes $\{X^{[V]}(t)\}$ by setting for every $V$
		\begin{equation}
		 X^{[V]}(t)=\frac{Y^{[V]}(t)}{V}.
		  \label{density}
		   \end{equation}
	\end{definition}
Let us remark that the name \emph{density process} originated from the population dynamics. In the reaction network model it represents the concentrations of the chemical species.
Following the theory of point processes \cite{AKbook,bremaud} density process $X^{[V]}(t)$ can be written in two equivalent (in a sense of the probability law) forms. The first form is the stochastic differential equation
\begin{equation}
dX_t=\sum_l \frac{l}{V} dM_l(t) \label{mc1}
\end{equation}
where $M_l(t)$ counts the occurrences of those reactions whose effect is to increase $Y^{[V]}(t)$ by $l$, and hence to increase the density process $X^{[V]}(t)$ by $l/V$. The state dependent rate associated with $M_l(t)$ is
\[ q^{[V]}_{X^{[V]}(t),X^{[V]}(t)+l/n}=V\; f_l \left(X^{[V]}(t)\right).\]

The second representation of the process $X^{[V]}(t)$  is obtained by substituting counting process $M_l(t)$ by independent unit-rate Poisson process $N_l(t)$. The effect of reactions with different speed is achieved by the \textit{time change}. The Poisson process implies a transformation that makes the individual time of each reaction to go `faster' when a higher jump rate is needed and `slower' otherwise. This leads to the following form of the process
\begin{align}
X_\ast^{[V]}(t) &= X_\ast^{[V]}(0)+\sum_{l} \frac{l}{V} N_{l}\left[ V \int_{0}^{t} f_{l}\left(X_\ast^{[V]}(s)\right) ds\right]\label{mc2}
%&=X^{[V]}(0)+ \sum_{l} \frac{l}{V} \left\{V \int_0^t f_l\left(X^{[V]}(s)\right) ds+ \tilde N_{l}\left[ V \int_{0}^{t} f_{l}\left(X^{[V]}(s)\right) ds\right]\right\}.\label{mc}
\end{align}
where $N_l(t)$ is an independent unit-rate Poisson process that counts the occurrences of the events which increase $Y^{[V]}(t)$ by $l$ (or the density process $X^{[V]}(t)$ by $\frac{l}{V}$). %The second line is a further rewriting which aimsWe denote the corresponding compensated processes by  $\tilde N_{l}(t)=N_{l}(t)-t$.

There are techniques \cite{lente,stewart1994}  to characterize both initial transient period and long run behaviour of the Continuous Time Markov Chains (CTMC). However, in practice if the state space of the CTMC is large, an analytical treatment is not feasible and an approximation is needed. The key idea is to construct a simpler process to approximate the original CTMC when it models the interaction of large groups. We briefly describe two of such approximations below.

\subsection{Approximations}
For large values of volume $V$, the jumps of the stochastic process \eqref{mc1} become more frequent and have a smaller magnitude suggesting that corresponding trajectories can be approximated by continuous functions (so-called the \emph{fluid} limit or the fluid approximations). In ~\cite{kurtz,kurtz1ode} a set of  \textit{ordinary} differential equations (ODE) providing the deterministic approximation of (\ref{mc1}) if both  volume and number of molecules are large is derived. This result is summarized below.

\begin{result}
\label{ode_th}
Let $x(t)$ be a deterministic solution of the $d$-dimensional ODE system
\begin{equation}\label{ode}
\dot x(t)=F(x(t))=\sum_{l\in C} l f_l(x(t))
\end{equation}
with initial condition $x(0)=x_0$. Let us assume that for each compact $K$ in the state space, the function $F$ is Lipshitz continuous in $K$ and that $\sum_l (|l|+|l|^2)\sup_{x\in K} f_l(x) <\infty$.
Let  $X^{[V]}(t)$ be as in \eqref{mc1} with the initial condition satisfying
\begin{equation}\lim_{V\rightarrow \infty}X^{[V]}(0) =x_0.\label{volume}\end{equation}
Fix time $T<\infty$. The density process $X^{[V]}(t)$ tends to $x(t)$ for all $t\leq T$ and
\begin{equation}
sup_{0\leq t \leq T}|X^{[V]}(t) - x(t)|= O \left(\frac{1}{\sqrt{V}} \right) 
\label{det_rate}
\end{equation}
with probability one as $V \to \infty$. The constant time horizon $T$ is arbitrary, but finite.
\end{result}

\noindent See  ~\cite{kurtz,kurtz1ode} for the proof. Let us remark that by equation \eqref{volume}, the parameter $V$ is related to the initial number of molecules and by letting $V$ to increase, the number of molecules in the system increases as well.
Since $x(t)$ provides a strong (path-wise) approximation of the process $X^{[V]}(t)$, for large $V$ every trajectory remains bounded in a small interval around the deterministic function $x(t)$. In such a regimen the stochastic nature of the process $X^{[V]}(t)$ is lost, with only the mean being relevant and approximated by $x(t)$ (note that for finite $V$ the mean of $X^{[V]}(t)$ is not necessarily given by $x(t)$, cf. \cite{jahnke2007solving}).
The limit (\ref{ode}) coincides with the classical deterministic formulation of the reaction network models (see e.g.   ~\cite{erditoth,feinberg,KurtzChemPhys}).

In a lot of cases (cf. \cite{lente}) the size of the system is not large enough to justify the deterministic approximation, and stochastic effects such as variance, skewness, bimodalities  are to be included into the approximating model.  A sharper continuous strong approximation that is able to capture stochastic fluctuations was obtained by  ~\cite{kurtz,Kurtz1976} in terms of the diffusion process. 

%\begin{result}
%\label{diffapprox}
%Let $X^{[V]}(t)$  be as in (\ref{mc1}). Let $G_\ast^{[V]}(t)$ be a diffusion process with initial condition satisfying $X^{[V]}(0)=G_\ast^{[V]}(0)$ and $\lim_{V\rightarrow \infty}X^{[V]}(0) =x_0$ which solves the following stochastic differential equation
%\begin{equation}
%G_\ast^{[V]}(t) = G_\ast^{[V]}(0)
%+\sum_{l} \frac{l}{V}\, \left[V\int_{0}^{t} f_{l}(G_\ast^{[V]}(s))ds+ W_{l}\left(V %\int_{0}^{t} f_{l}(G_\ast^{[V]}(s))ds\right)\right]\label{diff2}
%\end{equation}
%where $W_l(t)$ are independent standard Wiener process. \mozg{Let $\tau$ be the first time when process $G_\ast^{[V]}(t)$ hits the boundary.} Then, under some (rather general) conditions, it is possible to construct the processes $X^{[V]}_\ast(t)$ and $G_\ast^{[V]}(t)$ on the same probability space  such that  the difference between two paths fulfill
%\begin{equation}
%\sup_{0\leq t\leq \tau \wedge T} |X^{[V]}_\ast(t)-G_\ast^{[V]}(t)|= O \left(\frac{\log V}{V} \right) 
%\label{dif_rate}
%\end{equation}
% for any fixed time horizon $T$ with probability one as $V \to \infty$. 
%Note that the process $G_\ast^{[V]}(t)$ has the same law as the solution of the stochastic differential equation
%\begin{equation}
%G^{[V]}(t) = G^{[V]}(0)
%+\sum_{l} \frac{l}{\sqrt{V}}\, \left[\sqrt{V}\int_{0}^{t} f_{l}(G^{[V]}(s))ds+ \int_{0}^{t} \sqrt{f_{l}(G^{[V]}(s))} dW_{l}(s)\right]\label{diff1}
%\end{equation}
%due to the theory of time changed Wiener integrals (see, e.g., \cite{oksendal}, %Theorem 8.5.7).
%\end{result}

\begin{result}\label{diffapprox}
	Let $X^{[V]}(t)$  be as in (\ref{mc1}) and let $x(t)$ solve \eqref{ode} with initial condition $x(0)=x_0$. Let $G_\ast^{[V]}(t)$ be a diffusion process with initial condition satisfying $X^{[V]}(0)=G_\ast^{[V]}(0)$ and $\lim_{V\rightarrow \infty}X^{[V]}(0) =x_0$ which solves the following stochastic differential equation (given in the integral form)
	\begin{equation}
	G_\ast^{[V]}(t) = G_\ast^{[V]}(0)
	+\sum_{l} \frac{l}{V}\, \left[V\int_{0}^{t} f_{l}(G_\ast^{[V]}(s))ds+ W_{l}\left(V \int_{0}^{t} f_{l}(G_\ast^{[V]}(s))ds\right)\right]\label{diff2}
	\end{equation}
	where the $W_l(t)$ are independent standard Wiener processes.
	Let $E^{[V]}\subset \mathbb{R^d}$ be the smallest hyperrectangle (Cartesian product of $d$ intervals) that contains the discrete state space of $X^{[V]}(t)$.
	Let $U$ be any open connected subset of $E^{[V]}$ that contains $x(t)$ for every $0\leq t\leq T$. Let $\bar f_l=\sup_{x \in U} f_l(x)<\infty$ and suppose $\bar f_l=0$ except for finitely many $l$. Suppose $M>0$ satisfies both the two equations below for any $x,y\in U$
	\begin{equation}\begin{aligned}
	\abs{f_l(x)-f_l(y)}&\leq M\abs{x-y}\\
	\abs{F(x)-F(y)}&\leq M\abs{x-y}.
	\end{aligned}\label{Lip}\end{equation}
	Let $\tau_V=\inf\{t:X^{[V]}(t)\notin U \textup{ or }G_\ast^{[V]}(t)\notin U\}$. Note that $\mathbb{P}(\tau_V>T)\rightarrow 1$ for $V\rightarrow\infty$. Then for $V\rightarrow\infty$,
	
	\begin{equation}
	\sup_{0\leq t\leq \tau_V \wedge T} |X^{[V]}_\ast(t)-G_\ast^{[V]}(t)|= O \left(\frac{\log V}{V} \right) 
	\label{dif_rate}
	\end{equation}
for any fixed time horizon $T$.
\end{result}

\noindent See \cite{kurtz} for the proof and for a better estimate of the distance \eqref{dif_rate}. The statement of Approximation \ref{diffapprox} is quite complex, therefore, we provide some rephrasing of the main conclusion and the main assumptions.

Regarding the conclusion, Approximation \ref{diffapprox} states that it is possible to construct coupled trajectories of the two processes  $X^{[V]}_\ast(t)$ and $G_\ast^{[V]}(t)$ on the same probability space (using the same random numbers) in the way that the maximum distance between them is vanishing with a rate $\frac{\log V}{V}$ when  $V\rightarrow\infty$.

Regarding the assumption, some of them are technical, while others deserve to be discussed in more details. Firstly, the initial concentration is kept constant when the volume increases, the large systems with a huge number of molecules are approximated. Secondly, the assumptions on the functions $f_l(\cdot)$ are rather natural in the context of chemical kinetics as they prescribe that there is a finite number of reactions and none of them has an infinite speed. Finally, the approximation is only valid in any open set U that is contained in $E^{[V]}$ and that contains the whole trajectory of the deterministic approximation $x(t)$.
The introduction of such open set prevents both $X^{[V]}_\ast(t)$ and $G_\ast^{[V]}(t)$ from visiting the boundary of $E^{[V]}$.  The concentration of each chemical species can never become negative. In some example the concentration of a species is unbounded, in other it has an upper bound. In the case it may be unlimited, the introduction of such open set is needed since  the approximation will only work as far as the processes do not exceed any arbitrarily large but finite threshold (excluding explosions).
% Let us also remark that the Lipshitz conditions \eqref{Lip} may fail at infinity for many reasonable mass-action rates.
Moreover, in the case if the concentrations vanish the results of Approximation  \ref{diffapprox} would not hold any more. Let us remark that when $V$ is large enough both processes will be arbitrary close to $x(t)$ with high probability and visits of the boundaries will become less frequent (and absent in the limit). For the medium-large size systems visits to boundary might still be possible and the approximation would fail.  This is recognised as an important problem and has attracted a lot of attention in the literature \cite{angius2015approximate,ruth2017constrained,complex}.

Importantly, the process $G_\ast^{[V]}(t)$ has the same law as the solution of the stochastic differential equation
\begin{equation}
G^{[V]}(t) = G^{[V]}(0)
+\sum_{l} \frac{l}{\sqrt{V}}\, \left[\sqrt{V}\int_{0}^{t} f_{l}(G^{[V]}(s))ds+ \int_{0}^{t} \sqrt{f_{l}(G^{[V]}(s))} dW_{l}(s)\right]\label{diff1}
\end{equation}
due to the theory of time changed Wiener integrals (see, e.g., \cite{oksendal}, Theorem 8.5.7).

We would like to emphasize that both $x(t)$ and $G^{[V]}(t)$ provide the strong approximations of $X^{[V]}(t)$ and are not different in this sense.
The first term of (\ref{diff1}) is similar to the term in the deterministic approximation (\ref{ode}), but the second one adds noise and represents the stochastic nature of the process. The approximation $G^{[V]}(t)$ preserves a random behaviour of the process and corresponds to the lower rate of the error in (\ref{dif_rate}) compared to rate (\ref{det_rate}) for the deterministic fluid approximation.  As a result, this approximation can be applied in many cases where the deterministic one fails. A few examples are given in Section \ref{ex}.

The process (\ref{diff1}) is widely used to model chemical reactions (and well-known in Chemistry under the name of \textit{Langevin equations}, see \cite{gillespie}), mainly as a trick to speed up simulations. In our opinion, the diffusion approximation result obtained by \cite{kurtz} is not fully appreciated and deserves to be disseminated and applied more widely. Indeed, in addition to the guarantee that the laws of the processes $X^{[V]}(t)$ and $G^{[V]}(t)$ are similar, it gives the constructive procedure to generate discretized trajectories of the two processes $X^{[V]}_\ast(t)$ and $G_\ast^{[V]}(t)$ on the same probability space (i.e., with the same random numbers) that they stay close to each other \emph{trajectory by trajectory} with probability one. Since such construction is not given (to out best knowledge) explicitly in any work easily accessible to non-mathematicians, we provide it in the next section.

\section{Construction of paired trajectories of CTMC and diffusion approximation \label{sec:KMT}}

The constructions of $X^{[V]}_\ast(t)$ and $G_\ast^{[V]}(t)$ are built on two preliminary steps and one key argument. {Firstly, let $\tilde N(t)= N(t) -t$  be a compensated Poisson process with zero mean.} Note that $\tilde N(t)$ is a martingale and equation \eqref{mc2} can be written  as

\begin{align}
X_\ast^{[V]}(t) &=X_\ast^{[V]}(0)+ \sum_{l} \frac{l}{V} \left\{V \int_0^t f_l\left(X_\ast^{[V]}(s)\right) ds+ \tilde N_{l}\left[ V \int_{0}^{t} f_{l}\left(X_\ast^{[V]}(s)\right) ds\right]\right\}.\label{mc3}
\end{align}
\noindent Secondly, notice that the sole difference between equation \eqref{mc3} and equation \eqref{diff2} is that independent compensated Poisson process $\tilde N_l (t)$ is substituted by independent Wiener process $W_l(t)$. The key argument is a consequence of the \textit{KMT theorem}, named after the authors of \cite{KomlosI}. It states that paired trajectories of Wiener and  Poisson processes can be constructed on the same probability space such that the uniform distance between them is suitably controlled. Following \cite{kurtz} and \cite{KomlosI} we state the following Proposition. 
\begin{prop}\label{lemma:MBePP} 
	Given a Wiener process $W(t)$, a compensated Poisson process $\tilde N(t)$ can be constructed on the same probability space such that for any $\beta>0$ there exist positive constants $\lambda, \kappa$ and $c$  such that 
	\[\mathbb{P} \left(\sup_{t \leq \beta V} |\tilde N(t)-W(t)| \leq c \log V +x\right) \leq \kappa V^{-2}\text{e}^{-\lambda x}\]
for any $V>1$ and $x>0$.
\end{prop}

Given coupled trajectories of compensated Poisson process $\tilde N_l (t)$ and independent Wiener processe$W_l(t)$ constructed by Proposition \ref{lemma:MBePP}, it is a (non-trivial) technical matter to show that the uniform distance between $X_\ast^{[V]}(t) $ and $G_\ast^{[V]}(t)$ fulfils equation \eqref{dif_rate}. We start from the revisiting the construction needed to generate paired discretized sample paths of $\tilde N_l (t)$ and $W_l(t)$ and then we demonstrate how to build a discretization scheme for $X_\ast^{[V]}(t) $ and $G_\ast^{[V]}(t)$.

We would like to stress that a Poisson process can be seen as the partial sums of its increments and that the problem of approximating partial sums by Wiener process (strongly) has received a great attention in the literature. Strassen \cite{strassen1967} has used the Skorohod's embedding  scheme to provide the first construction. This construction, however, was shown to have not the best convergence rate \cite{csorgo1975}. Instead, the new construction based on the \textit{quantile transformation} of the increments of the original process was proposed by \cite{csorgo1975}. The quantile transformation of each value, however, was insufficient, while transforming blocks of increments proved a step in the right direction. The intuitive explanation is based on the central limit theorem which states that the sum of several independent and identically distributed random variables (under some conditions) tends to be normally distributed which makes the {quantile transformation} close to the identity. Therefore, it was proposed by \cite{csorgo1975}  to divide the values of process in blocks and to apply quantile transformations to sums in these blocks. The similar idea was used by \cite{KomlosI} and further extended to the quantile transformation into the individual blocks. The construction by \cite{KomlosI} was proved to achieve the best possible convergence rate and is provided below.

\subsection{Construction of paired Wiener and Poisson processes}

%\mozg{Historical story. Everything started from Strassen \cite{yyy}. Second. Third is Cirgo. BEST RATES EMERGE. It is important to mention that other constructions. The explanation behind becomes apparent once the construction is made so we give it after.

%Simple one to one transformation.

%Poisson converges in distribution to normal. If lambda is big then it is closer to normal. If the window is large. If window is small you accumulate the error. TWO SOURCES of ERROR. This is comprimise is optimal. 

%for example, a simple  quantile transformation of each point of the Wiener process, provide a greater convergence rate. Moreover, the KMT construction has the maximal possible convergence rate \cite{KomlosI}. The main reason is that the construction proceeds in blocks of the length $2^j$ and the method controls the error in each block.}

Importantly, the work by \cite{KomlosI} proves the \emph{existence} of coupled Poisson and Wiener processes and gives the construction of these processes. Precisely, given asequence of independent standard normal random variables $\{\bar W_i\}_{i=1\cdots N}$, it is possible to construct sequence of independent standard random variables $\{\bar N_i\}_{i=1\cdots k}$ with given distribution $F(x)$. It is also shown that the processes of the partial sums $T_n=\sum_{i=1}^n\bar{W}_i$ and $S_n=\sum_{i=1}^n\bar{N}_i$ fulfil 
\[\mathbb{P}\left(\sup_{1\leq n\leq k}|S_n-T_n|> C \log k +x\right)< K \text{e}^{-\lambda x}\]
for any arbitrary $x$, $n$ and for some positive constants $C$, $K$, $\lambda$ which depend on $F$ only.

As stated above, the KMT theorem by \cite{KomlosI} is constructive and gives an explicit algorithmic expression for the random variables $\{\bar N_i\}_{i=1\cdots k}$ in terms of the sequence $\{\bar W_i\}_{i=1\cdots N}$. This construction is also known as the  \emph{Hungarian construction}. Below we present its easily coded version allowing to simulate two discretized trajectories of  Wiener {process with drift} and Poisson process based on the same random numbers. Note that one can equivalently generate either (i) Wiener process with drift and Poisson process or (ii) Wiener process and compensated Poisson process. The goal of the representation below is pedagogical, thus we focus on the most straightforward implementation of the method rather than on computational costs or a memory usage. We  refer the reader to \cite{KomlosI} for the mathematical justifications.

Let us consider the time interval $[0,n \Delta]$ and its discretization with fixed step $\Delta$, $\{0,\Delta,2\Delta,\ldots,n\Delta\}$. We specialize the KMT Theorem to the case when the random variables $\{\bar N_i\}_{i=1\cdots n}$ are standardized Poisson increments 
\begin{align*}
\bar{N}_i = \frac{N_l(i\Delta)-N_l((i-1)\Delta) - \Delta}{\sqrt{\Delta}}.
\end{align*}
Then, Poisson process and Wiener process {with drift} having the same mean and variance can be obtained on the discretized time interval $[0, n\Delta]$ as
\begin{align}
N(k\Delta)  &=  \sum_{i=1}^k \left(\sqrt{\Delta}\bar{N}_i + \Delta \right), \quad \quad k=1, \dots, n,\label{eq:paired:P}\\
W (k\Delta) &= \sum_{i=1}^k \left( \sqrt{\Delta}\bar{W}_i + \Delta \right), \quad \quad k=1, \dots, n,\label{eq:paired:W}
\end{align}
where random variables $(\bar{W}_i)_{i=1}^n$ are distributed according to the standard normal distribution function with cumulative distribution function $\Phi$. The construction proceeds as follows. Given standardized Wiener increments $\{\bar{W}_1,\bar{W}_2,\ldots,\bar{W}_n\}$, we would like to find corresponding standardized Poisson increments $\{\bar{N}_1,\bar{N}_2,\ldots,\bar{N}_n\}$.

Without loss of generality, assume that the length of the trajectory $n$ can be written as $n=2^K$ where $K$ is positive integer.  Following the notation of \cite{KomlosI}, we {introduce the following quantities}
$$V_j=T_{2^j}, \ \ \ \ \ \ \ \ \ \ V_{j,k}=T_{(k+1)2^j}-T_{k2^j}, \ \ \ \ \ \ \ \ \ \ \tilde{V}_{q,k}=V_{q-1,2k}-V_{q-1,2k+1}.$$
As Wiener increments $\{\bar{W}_1,\dots,\bar{W}_n\}$ are already given, one can compute all of these quantities. The values of $V_{j,k}$ for all $j=0,1,\ldots, K-1$ and $k=1,\ldots,n-1$ can be written as elements of ${K \times (n-1)}$ dimensional matrix $\mathbb{V}$ with entries
\[
\begin{bmatrix}
    T_2-T_1 & T_3-T_2 & \ldots & \ldots & \ldots& \ldots& \ldots&T_{n}-T_{n-1} \\
    T_4-T_2 & T_6-T_4 & \ldots & \ldots & T_{n}-T_{n-2} & 0 & \ldots & 0 \\
    T_8-T_4 & T_{12}-T_8 & \ldots & T_n-T_{n-4} & 0   & \ldots & \ldots &  0  \\
    \vdots & \vdots & \vdots & \vdots & \vdots & \vdots & \vdots & \vdots \\
    T_n-T_{n-2^{K-1}} &0 & \ldots& \ldots& \ldots& \ldots& \ldots  & 0 \\
\end{bmatrix}_{K,n-1}
\]
Using the elements of $\mathbb{V}$, $\tilde{V}_{q,k}$ for all $q=1,2,\ldots, K-1$ and $k=1,\ldots,n-1$ can be found as elements of ${(K-1) \times (\frac{n}{2}-1)}$ dimensional matrix $\tilde{\mathbb{V}}$
\[
\begin{bmatrix}
   V_{0,2}-V_{0,3} & V_{0,4}-V_{0,5}& \ldots & \ldots  & \ldots & V_{0,n-2}-V_{0,n-1} \\
   V_{1,2}-V_{1,3} & V_{1,4}-V_{1,5}& \ldots & V_{1,\frac{n}{2}-2}-V_{1,\frac{n}{2}-1}  & \ldots & 0 \\
 \vdots & \vdots & \vdots& \vdots& \vdots  \\
 V_{K-2,2}-V_{K-2,3} & 0 & \ldots  & \ldots  &\ldots  & 0 \\
\end{bmatrix}_{K-1,\frac{n}{2}-1}
\]
The matrix $\tilde{\mathbb{V}}$ can be computed by Algorithm \ref{algg1}.

\begin{algorithm}[H]
    \caption{\label{algg1}Computing elements of matrix $\tilde{\mathbb{V}}$}
  \begin{algorithmic}[1]
    \FOR{$q=1,2,\ldots, K-1$}
      \STATE compute $V_{q,2k}$-$V_{q,2k+1}$ for all $k\leq \frac{n}{2^q}-1 $ and set elements of matrix with $k>\frac{n}{2^q}-1$ equal $0$
    \ENDFOR
    \medskip
  \end{algorithmic}
\end{algorithm}
Similarly, let us introduce the quantities
$$U_j=S_{2^j}, \ \ \ \ \ \ \ \ \ \ U_{j,k}=S_{(k+1)2^j}-S_{k2^j}, \ \ \ \ \ \ \ \ \ \ \tilde{U}_{q,k}=U_{q-1,2k}-U_{q-1,2k+1}.$$
\noindent {Note that $S_i$ are not yet known. In fact, the KMT computes $S_i$ using $U_i$ which are to be found using $V_i$.} Let us define matrices $\mathbb{U}$ and $\tilde{\mathbb{U}}$  similarly to $\mathbb{V}$ and $\tilde{\mathbb{V}}$ such that the entries of $\mathbb{U}$ and $\tilde{\mathbb{U}}$ have the same structure, respectively, but in terms of Poisson increments $\bar N_i$.

Since the goal of the method is to compute Poisson increments based on Wiener increments, we rephrase our goal by saying that we aim to compute the first line ($j=0$) of the matrix $\mathbb{U}$. Before Poisson increments can be computed the cumulative distributions function, conditional cumulative distributions function and corresponding quantile transformations should be defined as follows
$$\ F_j(x)=\mathbb{P} \left(U_j<x \right) \ \ \ \ \ \ \ \ \ \ \ \ \ \ \ \ \ \ \ \ \ \ \ \ \ \ \ \ \ \ F_q(x|y)=\mathbb{P} \left(\tilde{U}_{q,0}<x|U_{q_0}=y \right)  $$
$$G_j(t)= {\rm sup} \{x: F_j(x) \leq t \}  \ \ \ \ \ \ \ \ \ \ \ \ \ \ \ \ \ \ \ \ G_q(t|y)={\rm sup} \{x: F_q(x|y) \leq t \}  $$
Let us define $\rm{Poi}_\lambda(x)$ as distribution function of a Poisson r.v. with intensity $\lambda$. The cumulative distribution function $F_j(x)$ takes the form $$F_j(x)={\rm Poi}_{2^j \Delta} \left(\sqrt{\Delta}x+2^j\Delta\right).$$
%\mozg{Note that $U_j$ is the sum of $2^j$  Poisson r.v.s and is a Poisson random variable itself with rate $2^j \Delta$ and cumulative distribution function $F_j$ which is given explicitly.}

%For the Poisson process, $F_j(x)$ takes the form $F_j(x)={\rm Poi}_{2^j \Delta} \left(\sqrt{\Delta}x+2^j\Delta\right)$ as a sum of $2^j$ Poisson r.v.s with intensity $\lambda$ and distribution function ${\rm Poi}_\lambda$.
The conditional distribution function $ F_j(x | y)$ can be calculated observing that if $A$ and $B$ are independent Poisson random variables with intensity $2^{j-1} \Delta$, then
\begin{align*}
\mathbb{P} \left(  A-B < t \; \vert \; A+B= j \right) = 
\begin{cases}
0 & t < -j \\
\sum_{-j \leq i < t, \; i \leq j} \frac{ \mathbb{P}\left(  B = \frac{j-1}{2} , A =  \frac{j+i}{2} \right)   }{\mathbb{P}(A+B =j)}&  -j \leq t \leq j \\
1 & t > j.
\end{cases}
\end{align*}
Noticing that $\tilde{U}_{j,0}$ has the same distribution as $(A-B)/\sqrt{\Delta}$ and $U_j$ has the same distribution as $(A+B-2^j \Delta)/\sqrt{\Delta}$ leads to
\begin{align*}
F_j(x | y) =  \mathbb{P} \left(  A-B < \sqrt{\Delta} \; x \; \vert \; A+B = \sqrt{\Delta}\; y + 2^j \Delta\right).
\end{align*}
Then, the elements of matrix $\mathbb{U}$ are computed by Algorithm \ref{KMT_algorithm}

\begin{algorithm}[H]
\label{KMT_algorithm}
  \begin{algorithmic}[1]
         \STATE Compute $\bar N_1=G_0\left( \Phi \left( \bar W_1\right) \right)$
          \STATE Compute the first column of $\mathbb{U}$ using $U_{j,1}=G_j\left( \Phi \left( 2^{-\frac{j}{2}}V_{j,1} \right) \right)$
                     \FOR{$j,k$ such that $U_{j,k}$ is computed}
          \STATE Compute $\tilde{U}_{j,k}= G_j\left( \Phi \left( 2^{-\frac{j}{2}}\tilde{V}_{j,k} \right) | U_{j,k} \right)$
          \STATE Compute
          $$U_{j-1,2k}=\frac{1}{2} \left(U_{j,k}-\tilde{U}_{j,k} \right) $$
                    $$U_{j-1,2k+1}=\frac{1}{2} \left(U_{j,k}+\tilde{U}_{j,k} \right) $$
                    \ENDFOR{ \ when elements $U_{0,k}=N_{k+1}$ are found for all $k=1,\ldots,n-1$.}
    \medskip
  \end{algorithmic}
  \caption{KMT algorithm}
\end{algorithm}
\noindent Algorithm \ref{KMT_algorithm} computes the elements of matrix $\mathbb{U}$ from the last line ($j=K-1$) to the first one ($j=0$). While the equations are explicit, the order of elements computations might not be straightforward at the first glance. Therefore, we provide a pseudo-code for the computation in Algorithm \ref{A3}

\begin{algorithm}[H]\label{A3}
	\label{computations}
	\begin{algorithmic}[1]
		\STATE Set $c_1=c_2=1$
		\FOR{$u=1,2,\ldots,K-2$}
		\FOR{$v=1,\ldots,c_2$}
		\STATE Compute $U_{q-1,2c_1}$ and $U_{q,2c_1+1}$ for all $q=1,\ldots,K-u$
		\STATE Compute $\tilde{U}_{q-1,2c_1}$ and $\tilde{U}_{q-1,2c_1+1}$ for all $q=2,\ldots,K-u$
		\STATE $c_1=c_1+1$
		\ENDFOR{\ $v$}
		\STATE $c_2=2c_2$
		\ENDFOR{\ $u$ }
		\medskip
	\end{algorithmic}
	\caption{Computing elements of matrix ${\mathbb{U}}$}
\end{algorithm}

The processes needed to construct the original density dependent process $X_\ast^{[V]}(t)$ and the diffusion approximation $G_\ast^{[V]}(t)$ can be obtained by applying (\ref{eq:paired:P}) and (\ref{eq:paired:W}), respectively.
\subsubsection{Illustration}
To illustrate the construction we consider a toy example with $n=16$ ($K=4$) and $\Delta=1$. We simulate $16$ standard normal random variables $\left(\bar{W}_i\right)_{i=1}^{16}$ which are then truncated to
\begin{align*} [&-0.18, -0.93, -0.78, -1.65, -0.41, -1.10, -1.69,  2.52,  1.40,\\
  &0.18, -0.96,  1.26,  1.48,  0.52, -2.25,  0.47]^{\rm T}\end{align*}
for reproducibility.

Then, $\tilde{\mathbb{V}}$ takes the form presented in Table \ref{V} and $\mathbb{U}$ has the elements listed in Table \ref{U},

\begin{table}[ht]
\centering
\begin{tabular}{r|rrrrrrrrrrrrrrrr}
  \hline
 & $k=1$ & $k=2$ & $k=3$ & $k=4$ & $k=5$ & $k=6$ & $k=7$ \\
  \hline
  $q=1$ & 0.87 & 0.69 & -4.21 & 1.22 & -2.23 & 0.95 & -2.72 \\ 
  $q=2$ & -2.33 & 1.27 & 3.78  \\ 
  $q=3$ & 1.65  \\ 
   \hline
\end{tabular}
\caption{Elements of the matrix $\tilde{\mathbb{V}}$ for the illustrative example. The missing values are zeros.} \label{V} 
\end{table}

\begin{figure}[h!]
  \centering
    \includegraphics[width=0.75\textwidth]{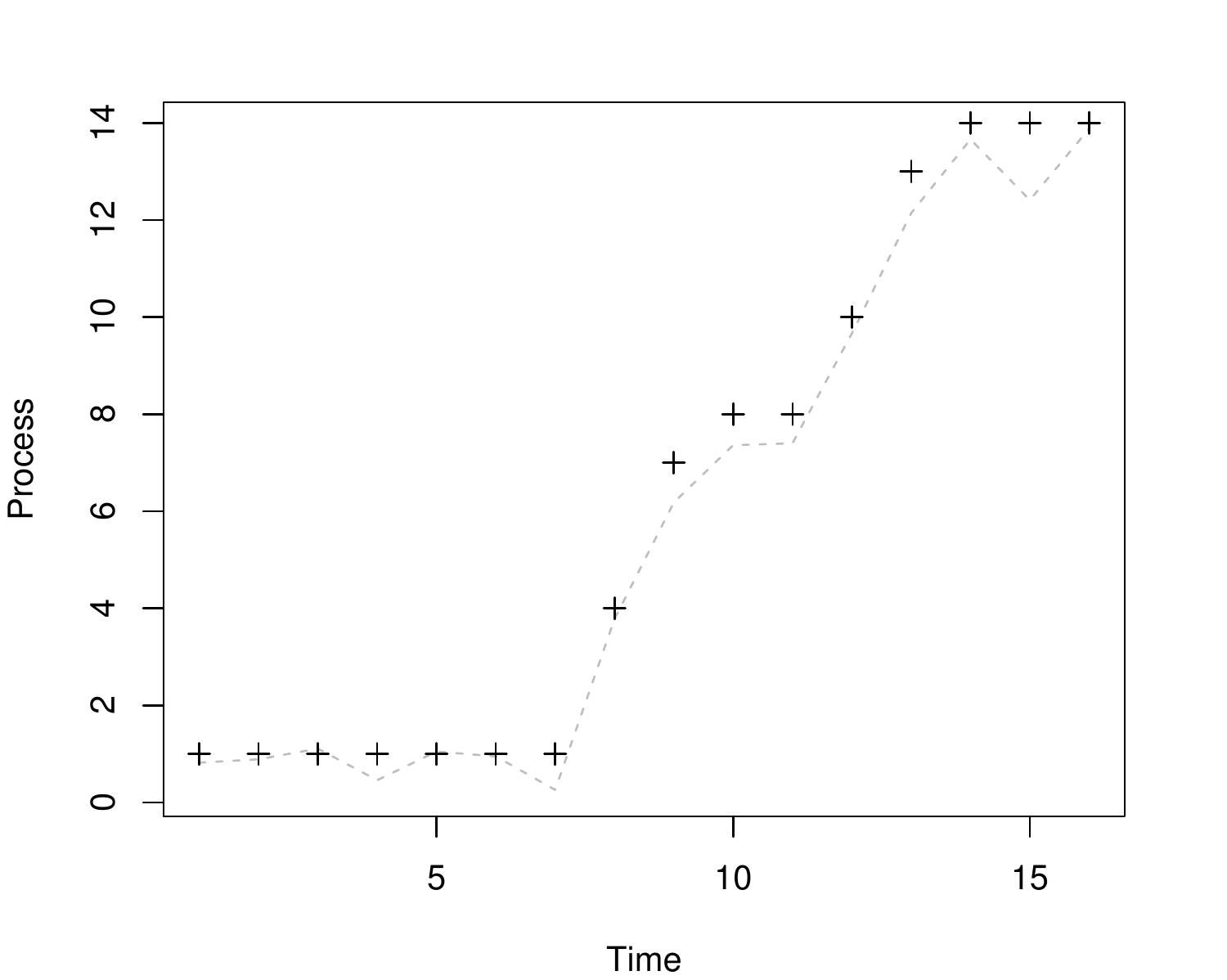}
    \caption{The pair of distretized Wiener process {with drift} (dashed grey line) $W_i$ and Poisson $N_i$ process (black cross)}
    \label{fig:illustration}
\end{figure}

%\begin{table}[ht]
%	\centering
%	\begin{tabular}{rrrrrrrrrrrrrrrr}
%		\hline
%		$j / k$ & $1$ & 2 & 3 & 4 & 5 & 6 & 7 & 8 & 9 & 10 & 11 & 12 & 13 & 14 & 15 \\ 
%		\hline
%		$0$  & \textbf{-1} & \textcolor{blue}{-1} & \textcolor{blue}{-1} & \textcolor{red}{-1} & \textcolor{red}{-1} & \textcolor{red}{-1} & \textcolor{red}{2} & 2 & 0 & -1 & 1 & 2 & 0 & -1 & -1 \\ 
%		$1$ & \textbf{-2} & \textcolor{blue}{-2} & \textcolor{blue}{1} & \textcolor{red}{2} & \textcolor{red}{0} & \textcolor{red}{2} & \textcolor{red}{-2} \\ 
%		$2$ & \textbf{-1} & \textcolor{blue}{2} & \textcolor{blue}{0} \\ 
%		$3$ & \textbf{2}  \\ 
%		\hline
%	\end{tabular}
%	\caption{Elements of the matrix $\mathbb{U}$ for the illustrative example. The missing values are zeros.} \label{U} 
%\end{table}
%\noindent where colors correspond to the order of computing. The procedure starts from the first column (black bold) and filling up next two columns (blue). Then each new column is used to obtain two more columns: columns 4-7 (red) and 8-15 (black), subsequently. It is easy to see that blocks of filling doubles (1, 2, 4 and 8 columns), respectively. This fact is coded using $c_2$ in Algorithm \ref{A3}.
\begin{table}[ht]
	\centering
	\begin{tabular}{rrrrrrrrrrrrrrrr}
		\hline
		$j / k$ & $1$ & 2 & 3 & 4 & 5 & 6 & 7 & 8 & 9 & 10 & 11 & 12 & 13 & 14 & 15 \\
		
		$0$  & \textbf{-1} & \underline{-1} & \underline{-1} & $\overline{-1}$ & $\overline{-1}$ & $\overline{-1}$ & $\overline{\ 2}$ & \textit{2} & \textit{0} & \textit{-1} & \textit{1} & \textit{2} & \textit{0} & \textit{-1} & \textit{-1} \\
		$1$ & \textbf{-2} & \underline{-2} & \underline{1} & $\overline{\ 2}$ & $\overline{\ 0}$ & $\overline{\ 2}$  & $\overline{-2}$  \\
		$2$ & \textbf{-1} & \underline{2} & \underline{0} \\
		$3$ & \textbf{2}  \\
		\hline
	\end{tabular}
	\caption{Elements of the matrix $\mathbb{U}$ for the illustrative example. The missing values are zeros. The bold, underlined, overlined and figures in italics correspond to the corresponding order of computing the elements.} \label{U}
\end{table}
\noindent where different fonts, under and over lines correspond to the order of computing. The procedure starts from the first column (black bold) and filling up next two columns (underlined). Then each new column is used to obtain two more columns: columns 4-7 (overlined) and 8-15 (in italics), subsequently. It is easy to see that blocks of filling doubles (1, 2, 4 and 8 columns), respectively. This fact is coded using $c_2$ in Algorithm \ref{A3}.
To compute the process of interest we apply (\ref{eq:paired:P}) and (\ref{eq:paired:W}). The obtained pair of processes is given in Figure \ref{fig:illustration}.
It is demonstrated that the constructed Poisson process captures the behaviour of the Wiener process with drift and the distance between these processes is controlled.

The explicit construction allows to provide a hint why the KMT construction has the best possible rate.  As the construction fits different blocks there are two kind of errors arise: (i) the sum of the errors of the quantile transformation and (ii) the maximum of the maxima of the partial sums in the individual block. Then, the construction chooses the optimal trade-off between these two errors. Furthermore, the simple quantile transformation $\hat{U}_{2n}$ of $\tilde{V}_{2n}$ is strictly independent of $U_{2n}$ and therefore the joint distribution of $\hat{U}_{2n}$ and ${U}_{2n}$ is not equal to the desirable one. This problem is solved by the conditional quantile transformation which fixes the value of $U_{2n}$. 

\subsection{Paired trajectories of the CTMC and of the diffusion approximation}
%Once the pair of processes $\{W_1,W_2,\ldots,W_n\}$ and $\{N_1,N_2,\ldots,N_n\}$ are constructed as shown above, 

As the processes $N_l(t)$ and $W_l(t)$ are continuous time processes, one can compute the discretised trajectories  $\hat{X}^{[V]}(t)$ of density process $X_\ast^{[V]}(t)$ given in (\ref{mc2}) using an Euler scheme with step $\delta$ (which might not to coincide with step $\Delta$ in Section 3.1). We would get
\begin{align}\label{eq:paired:Xdiscrete}
\hat{X}^{[V]} (j \delta)  = & \; \hat{X}^{[V]}  ((j-1) \delta) + \\
& + \sum_l \; \frac{l}{n} \left[  N_l \left(  n\delta \sum_{k=0}^{j-1} f_l(\hat{X}^{[V]} (k\delta)) \right) - N_l \left(  n\delta \sum_{k=0}^{j-2} f_l(\hat{X}^{[V]} (k\delta)) \right)  \right], \nonumber
\end{align}
for $j=1,\ldots,N$ with $\hat{X}^{[V]}(0)={X}^{[V]}(0)$. Similarly, one can obtain the discretised trajectories $\hat{G}^{[V]}(t)$ of the diffusion approximation ${G_\ast}^{[V]}(t)$ given in \eqref{diff2} by
\begin{align}\label{eq:paired:Zdiscrete}
\hat{G}^{[V]} (j \delta)  = & \; \hat{G}^{[V]} ((j-1) \delta) + \nonumber\\
& + \sum_l  \frac{l}{n} \;  \left[  W_l \left(  n\delta \sum_{k=0}^{j-1} f_l(\hat{G}^{[V]}(k\delta)) \right) - W_l \left(  n\delta \sum_{k=0}^{j-2} f_l(\hat{G}^{[V]}(k\delta)) \right)  \right] ,
\end{align}
for $j=1,\ldots,N$ with $\hat{G}^{[V]}(0)={G}_\ast^{[V]}(0)$. Since processes $N_l(t)$ and $W_l(t)$ are not available in continuous time, but only on a discrete grid of amplitude $\Delta$, one needs to introduce a further approximation by replacing the four times
\[\begin{aligned}
 {}&n\delta \sum_{k=0}^{j-1} f_l(\hat{X}^{[V]} (k\delta)), \quad  n\delta \sum_{k=0}^{j-2} f_l(\hat{X}^{[V]} (k\delta))\\
 &n\delta \sum_{k=0}^{j-1} f_l(\hat{G}^{[V]}(k\delta)), \quad n\delta \sum_{k=0}^{j-2} f_l(\hat{G}^{[V]}(k\delta))
\end{aligned}\]
by the closest times on the grid  of obtained Wiener and Poisson processes. This would often require long trajectories of $N_l$ and $W_l$ using extremely small step $\Delta$ to get $\hat{G}^{[V]}(t)$, $\hat{X}^{[V]}(t)$ trajectories of a moderate length. However, the challenge is computational only and can be resolved by storing long trajectories of these processes. As a final remark, we would like to emphasize that we further consider two reaction network examples for which the trajectories of CTMC and its strong diffusion approximation are provided (see Figure \ref{fig:3} and Figure \ref{fig:4}). However, due to the computation costs, it is strongly recommended to a reader to simulate independent trajectories of \eqref{mc1}  and \eqref{diff1} using the classical algorithms instead if the trajectory-by-trajectory behaviour is not of interest.

\section{Examples}\label{ex}
In this section we describe two examples of the chemical reaction systems taken from the recent literature. The aim is to discuss the ability of the deterministic and diffusion approximations to capture the dynamical properties of the original Markov Chain.
\subsection{A toy model of metabolism and an interpretation of the deficiency} \label{ex1}
The deficiency of the network has been introduced as an algebraic property of the reaction graph in Section \ref{prima}. The authors of \cite{polettini2015} proposed a thermodynamic interpretation of the deficiency in terms of the entropy balance.
According to this interpretation, the deficiency can be understood as a number of the `hidden' closed pathways, or thermodynamic cycles. In case $\theta=0$, the average stochastic dissipation rate equals the rate of the corresponding deterministic model. They proposed the following toy model inspired by metabolism for the illustration
\begin{equation} \label{system}\begin{gathered}
\mathrm{nE} \xrightleftharpoons[\lambda_6]{\lambda_5} \mathrm{\oldemptyset} \xrightleftharpoons[\lambda_2]{\lambda_1} \mathrm{N} \\
\mathrm{N+mE} \xrightleftharpoons[\lambda_4]{\lambda_3} \mathrm{(m+n)E} 
\end{gathered} \end{equation}
\noindent where $\mathrm{N}$ is number of nutrients and $\mathrm{E}$ is number of tokens of energy. The first reaction introduces (eliminates) nutrients and energy to (from) the environment. The second reaction processes the nutrients and $m$ tokens of energy to produce more energy and vice versa. Following \cite{polettini2015}, we fix $n=2$. The stoichiometric matrix 
\begin{equation} \label{sm}
\begin{bmatrix}
1& -1 & -1 & 1 & 0 & 0 \\
0 & 0 & 2 & -2 & -2 & 2 \\
\end{bmatrix}
\end{equation}
displays in the $i$-th column the increment  caused in $(\mathrm{N},\mathrm{E})$ by the reactions with propensities $\lambda_1,\ldots,\lambda_6$ in system \eqref{system}. The approximate rates of reactions (neglecting the terms with higher order in $1/V$ in equation \eqref{massaction}) equal 
$$q^{(1)}_{(N,E),(N+1,E)}=\lambda_1 V, \ \ \ \ \  q^{(2)}_{(N,E),(N-1,E)}=\lambda_2 N $$
$$q^{(3) }_{(N,E),(N-1,E+2)}=\lambda_3 \frac{N E^m}{V^m}, \ \ \ \ \ q^{(4)}_{(N,E),(N+1,E-2)}=\lambda_4 \frac{E^{2+m}}{V^{1+m}}.$$
$$q^{(5)}_{(N,E),(N,E-2)}=\lambda_5 \frac{E^2}{V}, \ \ \ \ \  q^{(6)}_{(N,E),(N,E+2)}=\lambda_6 V$$

If $m$ is strictly positive, the network is made of 5 complexes with 2 connected components and the stoichiometric space has a dimension of 2. Then, the deficiency equals $\theta=5-2-2=1$ and is non-vanishing. In contrast, if $m=0$, the network is made of just 3 complexes, it has the single connected component and the stoichiometric space has a dimension of 2. Thus, there is no  deficiency in the system $\theta=3-1-2=0$.

\begin{figure}[h!]
	\centering
	\includegraphics[width=1\textwidth]{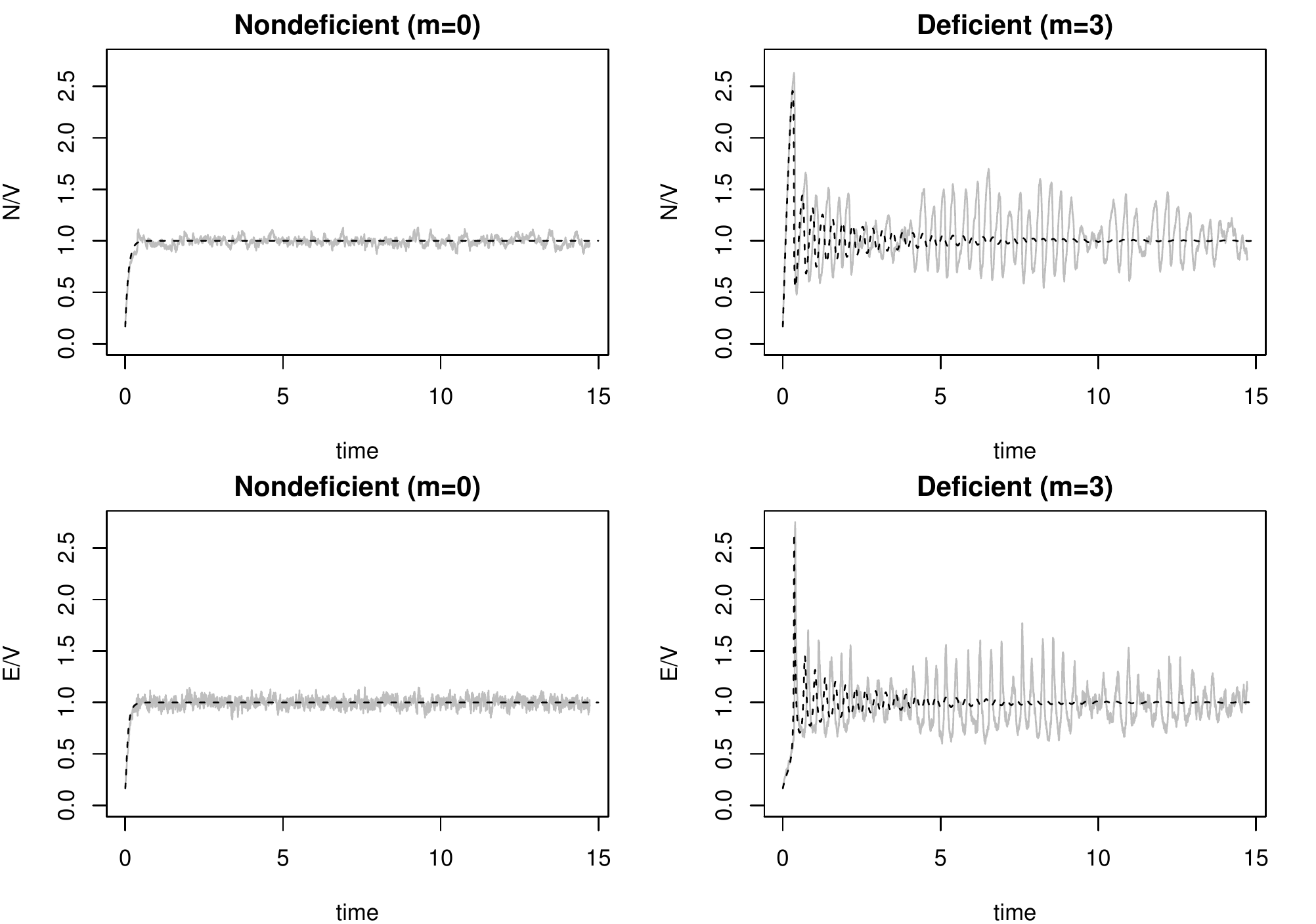}
	\caption{Deterministic (dashed black) and stochastic (solid grey) trajectories of N/V (upper part) and E/V (lower part) in the system (\ref{system}) in the nondeficient case (left panels, m=0) and in the deficient case (right panels, m=3).
%The two models disagree in the deficient case. According to the deterministic model, oscillations are damped, while they are sustained in the stochastic system.
	 }
	\label{fig:1}
\end{figure}

Following the choice of parameters by \cite{polettini2015}, we set $\lambda_1=10, \lambda_2=1, \lambda_3=10, \lambda_4=1, \lambda_5=10, \lambda_6=1$. Trajectories of the stochastic model (CTMC) of the system \eqref{system} in both non-deficient ($m=0$) and deficient ($m=3$) cases are given in Figure \ref{fig:1}. Figure \ref{fig:1} also shows the deterministic approximation which solves the system of ODEs
$$\dot{u}=10-u-10ue^m+e^{n+m} $$ 
$$\dot{e}=2\left(10ue^m-e^{n+m} + 1 - 10e^n \right) $$
where the variable $(u,e)$ are interpreted as nutrients concentrations $u=N/V$ and energy concentration $e=E/V$. {We fix the value $V=600$ as in the original example}. To generate the trajectory of the CTMC we use the stochastic simulation algorithm by \cite{Gillespie_1977} implemented in \texttt{R} \cite{Rurl}.

While being quite accurate in the non-deficient case, the deterministic approximation fails to catch properties of the stochastic system in the deficient case {for the chosen value of $V$}. Indeed, according to the deterministic model, for $m=3$ the system should display damped oscillation around the equilibrium that becomes of negligible amplitude as time goes. At the same time, the stochastic model prescribes sustained oscillations that are reducing their amplitude. This important qualitative feature is missed by the deterministic model. {We remark that the inadequacy of the deterministic approximation is due to the choice of $V$. While the result in Approximation \ref{ode_th} says that the deterministic approximation is valid for $V$ large enough, it, however, fails to reflect the properties of the original process for the original choice of $V$.}

\begin{figure}[h!]
	\centering
	\includegraphics[width=1\textwidth]{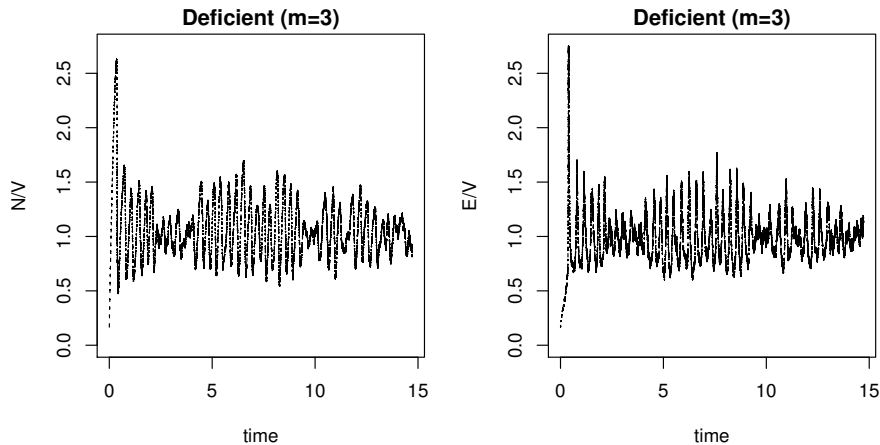}
	\caption{Trajectories of N/V (left) and E/V (right) according to the diffusion approximation \eqref{diff1}, for the system (\ref{system}) only in the deficient case  m=3. The diffusion approximation is in very good agreement with the Markov Chain of Figure \ref{fig:1}.}
	\label{fig:2}
\end{figure}

We further investigate whether the diffusion approximation is able to capture the qualitative dynamical properties of the system. Trajectories of the diffusion approximation \eqref{diff1} are given in Figure \ref{fig:2}.

\begin{figure}[h!]
	\centering
	\includegraphics[width=1\textwidth]{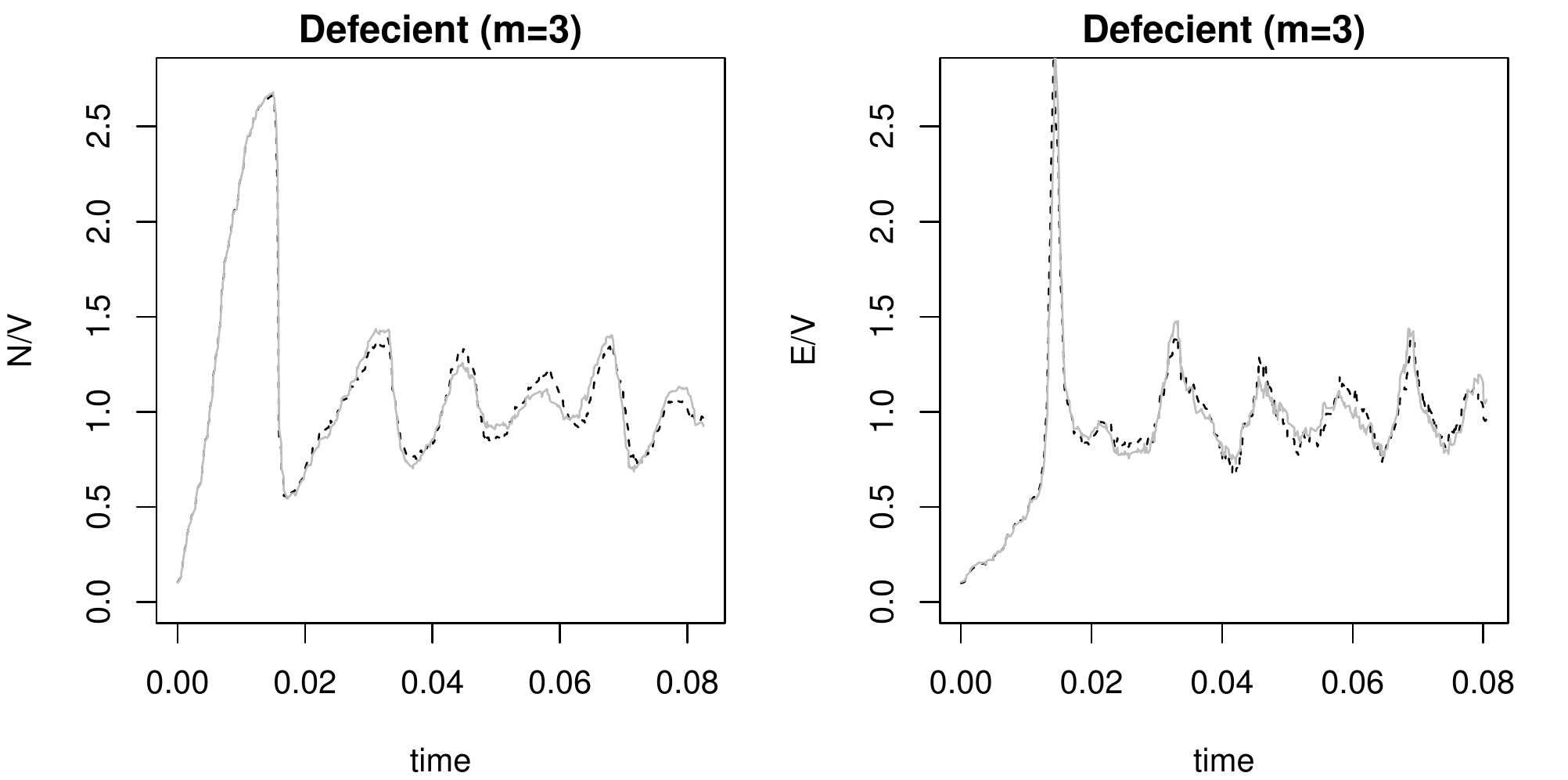}
	\caption{Paired trajectories of N/V (left)  and E/V (right)  according to the CTMC, (grey line) and the diffusion approximation, (dashed black line) for system (\ref{system}) in the deficient case $m=3$.
%The diffusion approximation is in very good agreement with the Markov Chain. For computational reasons time is limited to 2.
}
	\label{fig:3}
\end{figure}

One can see that oscillations are not damped according to the diffusion approximation. The trajectories of the diffusion closely resemble the behaviour of the original Markov chain. Importantly, the computation cost to compute the diffusion approximation is significantly lower than for the original CTMC. {For instance, to obtain a trajectory up to $t=50$, CTMC takes nearly 14.7 seconds, while the diffusion approximation takes less than $1.8$ seconds.} Note, however, that trajectories in Figure \ref{fig:2} are generated independently of CTMC, by Euler-Maruyama discretization method, applied to the equation \eqref{diff1}. To check how well the diffusion can approximate the trajectory of the CTMC given the same random numbers, we apply the algorithms described in Section \ref{sec:KMT} to generate paired trajectories. Due to the high computational costs, we  limit the time to 2 which would be enough to see the general pattern. Two paired trajectories are given in Figure \ref{fig:3}.

One can see that the corresponding trajectories are located close to each other to a very high extent and are in agreement. It follows that the diffusion approximation can mimic the behaviour of the original process, but with much less computational and analytical costs. Finally, we would like to outline that the diffusion approximation should not be considered as a short cut (to reduce the simulation time) only, as the theoretical analysis of systems with oscillations can been also performed on its basis \cite{baxendale2011}.  

\subsection{A \emph{minimal} chemical reaction systems with bistability} \label{ex2}
A bistable system is a system which has two stable equilibrium states and can be resting in either of these states. Bistable systems have been studied extensively to analyse kinetics, non-equilibrium thermodynamics and stochastic resonance. Due to its outstanding importance, the theoretical foundations of the bistability such as necessary and sufficient conditions have attracted an extensive attention in the literature, see, e.g., \cite{joshi2013atoms,wilhelm2009}.

\begin{figure}[h!]
	\centering
		\includegraphics[width=1\textwidth]{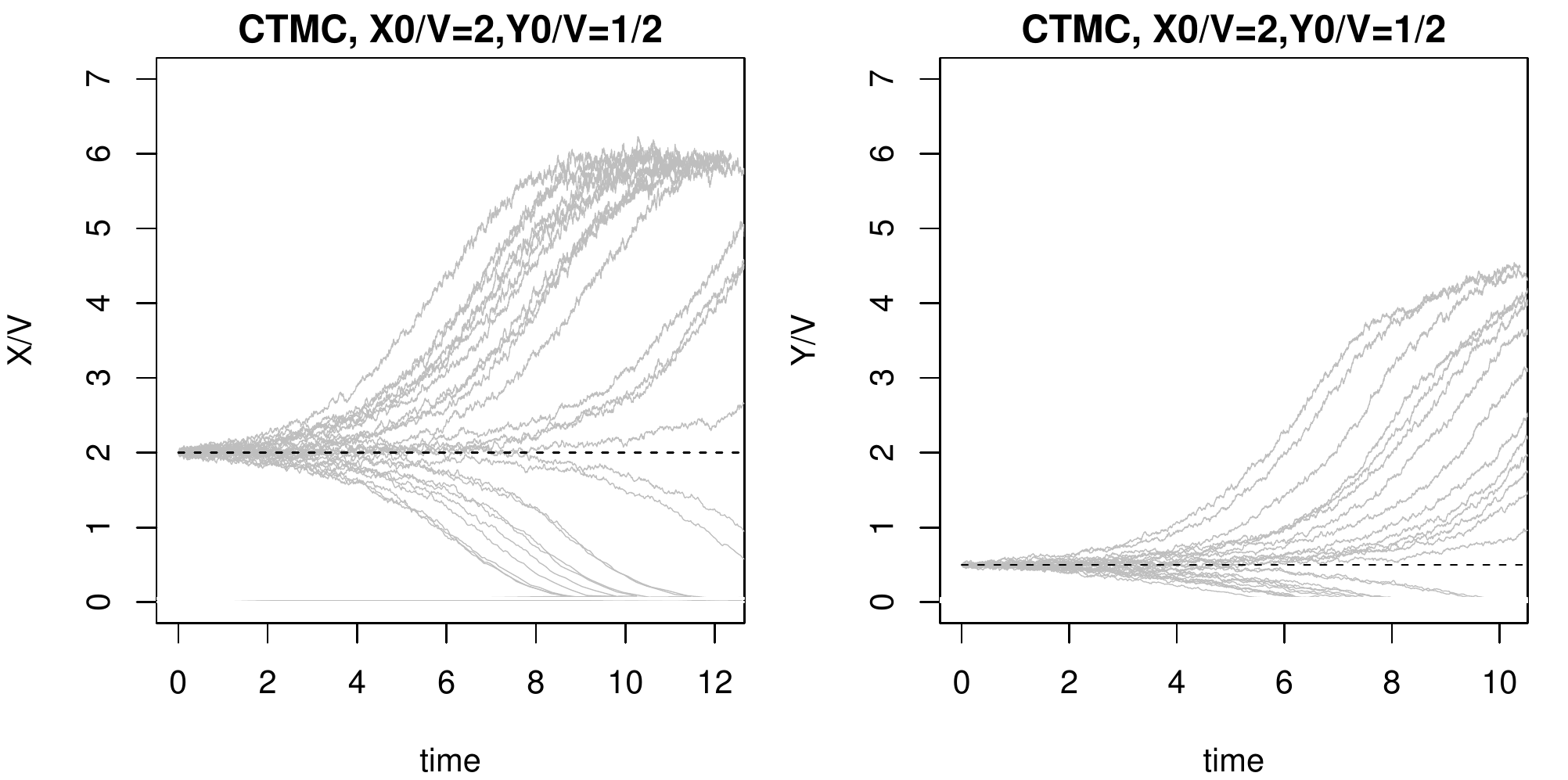}
	\caption{100 stochastic (grey line) and 1 deterministic (black dashed line) trajectories  of $X/V$ (left panel) and $Y/V$ (right panel) in system (\ref{bistable_system}) starting at the unstable steady state $x(0)=2, y(0)=1/2$.}
	\label{fig:bistable}
\end{figure}

The approach to formulate the necessary conditions of the bistability proposed by \cite{wilhelm2009} is to find a corresponding minimal bistable chemical system (MBCS). The authors use the wording \emph{chemical} system to indicate a special case of a \emph{mass-action} system such that all the reaction involved are at most bimolecular. More complicated reactions are indeed believed not to be \emph{physical}. We refer a reader to the original proposal of \cite{wilhelm2009} for the detailed definition of the minimal chemical system and for the comparison to alternative definitions, for instance, Schlogl model, \cite{schlogl1972}. The proposed MBCS consists of four reactions
\begin{gather}
\label{bistable_system}
\mathrm{S}+\mathrm{Y} \stackrel{\lambda_1}{\longrightarrow} \mathrm{2X} \nonumber \\
\mathrm{2X} \stackrel{\lambda_2}{\longrightarrow}  \mathrm{X+Y} \\
\mathrm{X+Y} \stackrel{\lambda_2}{\longrightarrow}  \mathrm{Y+P} \nonumber \\
\mathrm{X} \stackrel{\lambda_4}{\longrightarrow}  \mathrm{P} \nonumber
\end{gather}
where $X$, $Y$ are reactants and $S$, $P$ are substrates and products whose concentrations are kept fixed. The corresponding stoichiometric matrix takes the form
\[
\begin{bmatrix}
2 & -1 & -1 & -1 \\
-1 & 1 &  0 &  0 \\
\end{bmatrix}
\]
with the same convention adopted for equation \eqref{sm}
and the approximate rates (neglecting the terms with higher order in $1/V$ in equation \eqref{massaction}) are
$$q^{(1)}_{(X,Y),(X+2,Y-1)}=\lambda_1 Y, \ \ \ \ \  q^{(2)}_{(X,Y),(X-1,Y+1)}=\lambda_2 \frac{X^2}{V} $$
$$q^{(3)}_{(X,Y),(X-1,Y)}=\lambda_3 \frac{XY}{V}, \ \ \ \ \  q^{(4)}_{(X,Y),(X-1,Y)}=\lambda_4 X$$
where the constant concentration of $S$ is incorporated in $\lambda_1$. The system can be described by the ODE system
$$\dot{x}=2\lambda_1 y - \lambda_2 x^2 - \lambda_3 x y - \lambda_4 x $$
$$\dot{y} = \lambda_2 x^2 - \lambda_1 y $$
where $x$ and $y$ are concentrations of $X$ and $Y$, respectively. Setting $\lambda_2=1$, without restrictions of generality \cite{wilhelm2009} has shown that the system has three steady states $\bar{x}_1=\bar{y}_1=0$ and $\bar{x}_{2,3}=\frac{\lambda_1 \pm \sqrt{\lambda_1 D}}{2\lambda_3}, \ \bar{y}_{2,3}=\frac{\bar{x}_{2,3}^2}{\lambda_1}$ 
where $D= \lambda_1 - 4 \lambda_3 \lambda_4$ and the steady states 1 and 3 are stable and the steady state 2 is unstable. Following \cite{wilhelm2009} we set $\lambda_1=8$, $\lambda_2=\lambda_3=1$ and $\lambda_4=1.5$ for the illustration. In this case the steady states are $\bar{x}_1=\bar{y}_1=0$, $\bar{x}_{2}=2$, $\bar{y}_{2}=1/2$  and  $\bar{x}_{3}=6$, $\bar{y}_{3}=9/2$. One hundred trajectories of the CTMC corresponding to system (\ref{bistable_system}) and starting at the unstable steady state $x(0)=2, y(0)=1/2$ are given in Figure \ref{fig:bistable} where the deterministic approximation of the system is also presented.

The different trajectories originating at the same unstable steady state are driven by the noise, towards one of the stable equilibria picked randomly. Let us remark that despite many trajectories are initially attracted to the upper equilibrium, sooner or later they will escape its domain of attraction and they will end up visiting the state $(0,0)$ that is absorbing. Notice that this effect cannot be illustrated in the simulations since the time required to leave the upper equilibrium is much larger than the time windows that one can explore.

\begin{figure}[h!]
	\centering
		\includegraphics[width=1\textwidth]{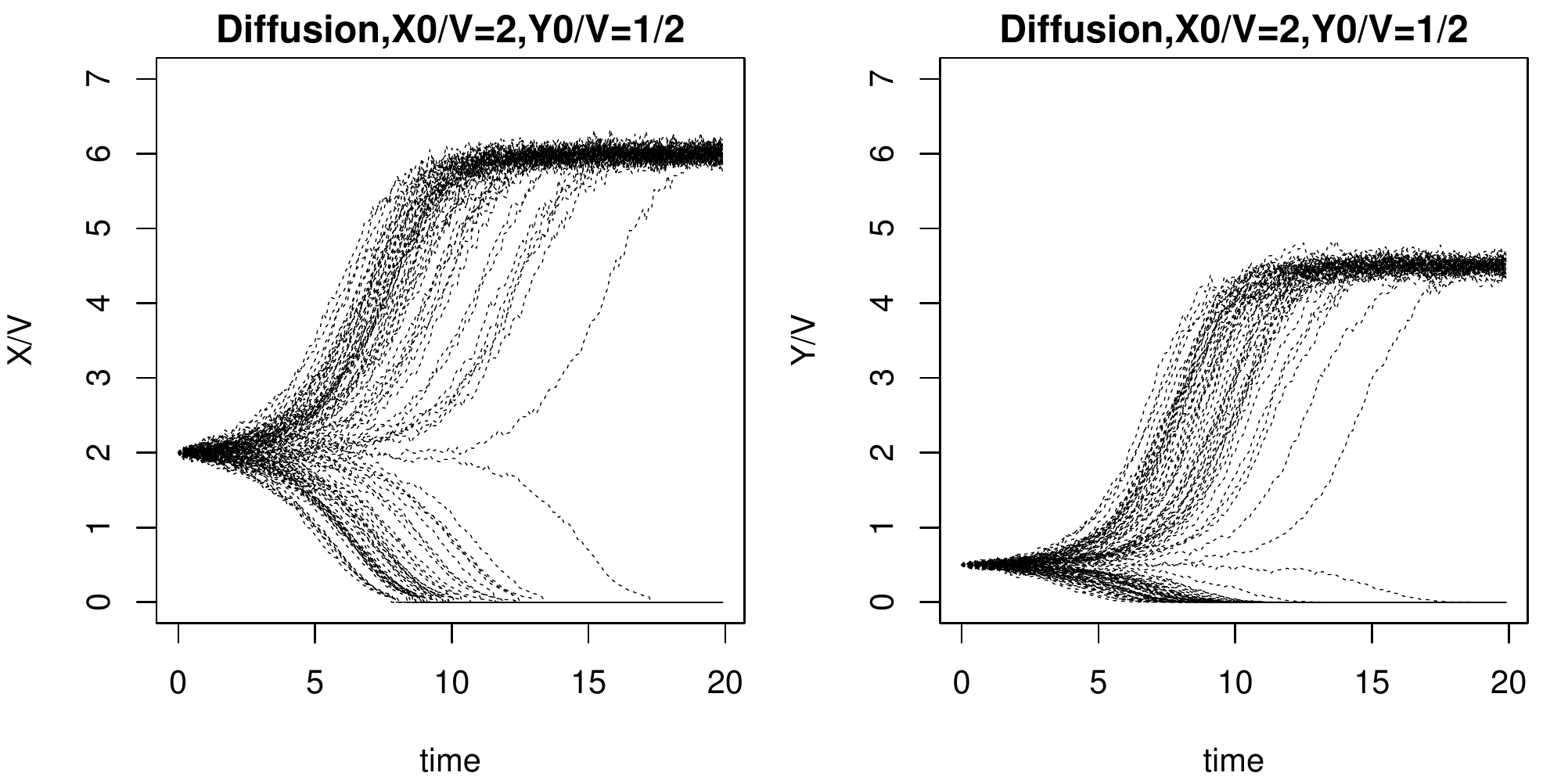}
	\caption{100 discretized trajectories (dashed black lines) of the diffusion approximation of $X/V$ (left panel) and $Y/V$ (right panel) starting at the unstable steady state $x(0)=2, y(0)=1/2$.}
	\label{fig:approx_bistable}
\end{figure}

 The deterministic approximation is not able to capture this complex and rich behaviour of the system.  Therefore, this approximation could lose important properties of the original process and should not be used. We then investigate the behaviour of the diffusion approximation. One hundred discretized trajectories (starting at the same point of unstable steady state) of the diffusion approximation are given in Figure \ref{fig:approx_bistable}.

The diffusion approximation mimics the qualitative  behaviour of the original CTMC very closely. Further, we study whether the diffusion approximation is able to reproduce the behaviour of the original process for different starting points. We calculate the proportion of trajectories of the CTMC and of the diffusion approximation that is attracted to each steady state for different initial conditions at some fixed time point $t$.  The results of nine sets of the initial points and the fixed time $t=20$ for the CTMC and the diffusion approximation are given in Table \ref{tab:bistable}.
\begin{table}[ht]
	\centering
	\caption{\label{tab:bistable} The proportion of times the CTMC (upper lines) and the diffusion approximation (lower lines) is attracted to the first steady state. The results are based on $10^4$ replications.} 
	\begin{tabular}{cccc}
		\hline
		Initial point & $y(0)=0.45$ & $y(0)=0.5$ & $y(0)=0.55$ \\
		\multirow{2}{*}{$x(0)=1.95$}  & 94.89\% & 75.24\% & 34.85\% \\
		& 94.89\% &  75.54\% & 34.49\% \\
		\\
		\multirow{2}{*}{$x(0)=2.00$}  &  85.12\% &  49.63\% & 14.76\% \\
		& 85.23\% &  49.55\% & 14.48\% \\
		\\
		\multirow{2}{*}{$x(0)=2.05$}  & 65.09\% &  25.04\% & 4.83\% \\
		& 65.52\% & 25.45\% & 4.63\% \\
		\hline
	\end{tabular}
\end{table}
Clearly, the diffusion approximation correctly reflects the behaviours of the original process up to the moment at which the absorbing state $(0,0)$ (a boundary of the state space) is reached, as described in the comments following the statement of Approximation 2. Importantly, the computation time for $100$ trajectories of the CTMC was nearly $96$ minutes, while the diffusion approximation took a half of the minute. Both this and its good property to mimic the behaviour of the original process make the diffusion approximation a reasonable tool to study the behaviour of the minimal bistable chemical system.

\begin{figure}[h!]
	\centering
		\includegraphics[width=1\textwidth]{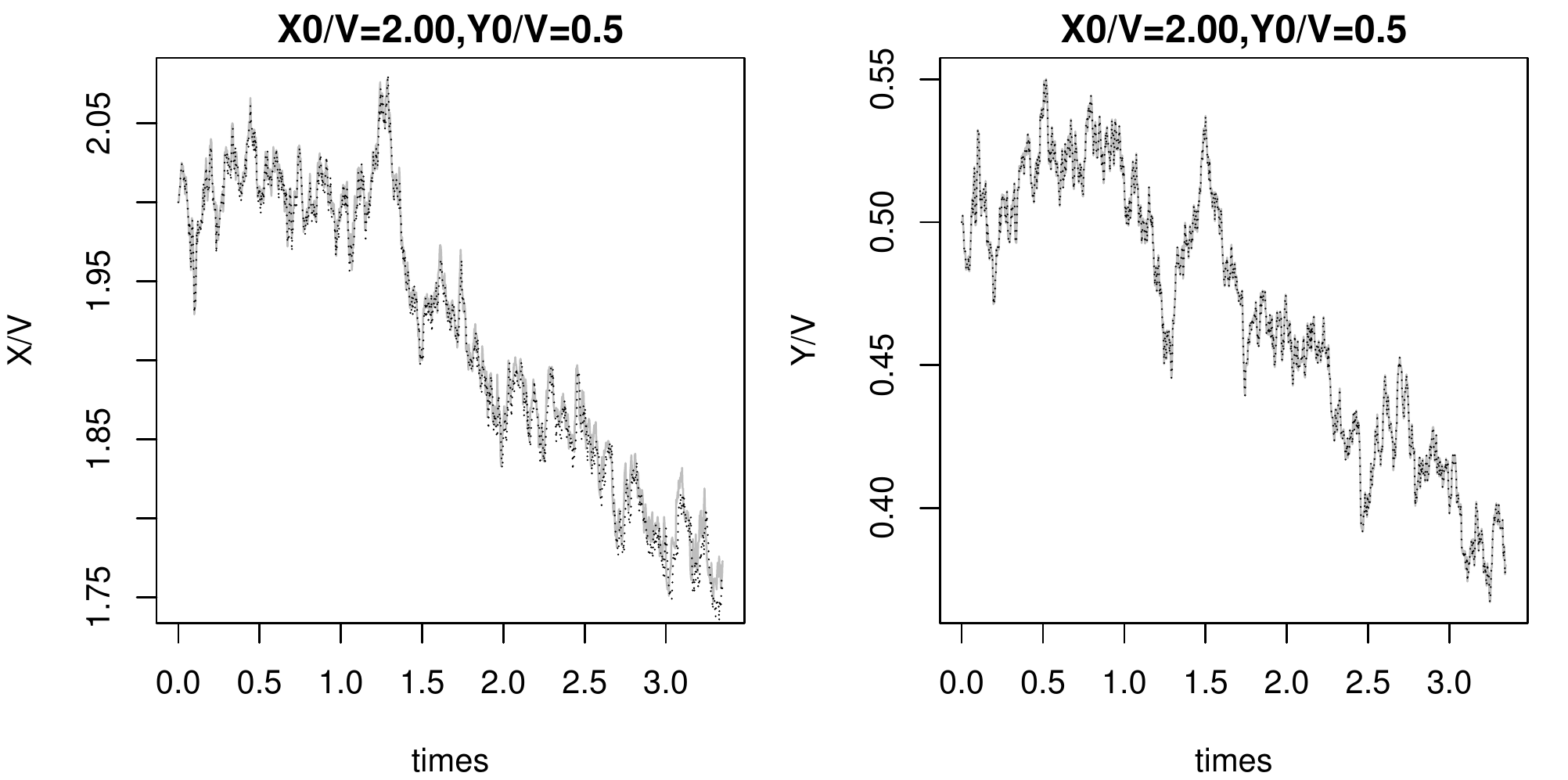}
	\caption{Paired discretized trajectories of the diffusion approximation (dashed black) and of the original CTMC (grey solid) of $X/V$ (left panel) and $Y/V$ (right panel) starting at the unstable steady state $x(0)=2, y(0)=1/2$.}
	\label{fig:4}
\end{figure}

As a further investigation on the diffusion approximation, we now compare the paired trajectories of the original process and the diffusion. Again, we limit time to $t=3.5$ due to the computational cost. The two trajectories are plotted in Figure \ref{fig:4}.
One can see that the both trajectories show a great agreement on the whole trajectory. The distance between processes stays little. The same conclusions were obtained for different starting values and, therefore, are not provided here.

\section{Limitations and perspective \label{discuss}}
In this work we demonstrated that the deterministic and diffusion approximations are useful tools in the modelling of reaction networks. The diffusion approximation is able to capture the behaviour of the original process and is able to mimic  trajectories of the CTMC for many different system. However, many questions of their applicability to important problems remain unanswered. Firstly, both approximations are derived on a finite time horizon. It is well known that deterministic equations may fail to catch the limiting distribution of the corresponding stochastic model when time goes to infinity as it happens in the Example presented in Section \ref{ex2}, and for all the chemical systems with absolute concentration robustness \cite{acr,danieleacm}. At the same time, they can capture such asymptotic behaviour correctly in case of the complex balanced stochastic systems \cite{danielecarsten}.  Similar results are not yet obtained for the diffusion approximation. Secondly, both the diffusion and the deterministic approximations are known to fail when the state space of the processes is bounded and the boundaries are visited with non-negligible probability. This may be a major drawback for the medium-large size systems where the size is not large enough. Alternative approximations have been proposed for this case in \cite{angius2015approximate,ATPN14,bbs,ruth2017constrained,complex}, but the complete mathematical theory is still under development.

%\begin{acknowledgements}
%If you'd like to thank anyone, place your comments here
%and remove the percent signs.
%\end{acknowledgements}

\bibliographystyle{spmpsci}
\bibliography{bibl}  

\end{document}